\documentclass[10pt,a4paper]{article}
\usepackage{amsmath}
\usepackage{amsthm}
\usepackage{amsfonts}
\usepackage{amssymb}
\usepackage[utf8]{inputenc}
\usepackage{hyperref} 
\usepackage{graphicx}

\theoremstyle{plain}
\newtheorem{theorem}{Theorem}

\newtheorem{claim}{Claim}[theorem]

\theoremstyle{definition}

\newtheorem{conjecture}[theorem]{Conjecture}

\begin{document}

\title{Note: a counterexample to a conjecture of Jackson about hamiltonicity of diregular digraphs.}

\author{Georgi Guninski\thanks{email :{\texttt{guninski@guninski.com}}}
}

\maketitle

\begin{abstract}
$3$-diregular circulant digraph on $12$ vertices is a counterexample of Jackson's conjecture about hamiltonicity of diregular digraphs
\end{abstract}

\begin{conjecture}[Jackson \cite{hand1} p. 68, 5.42]
Every $k$-diregular oriented graph on at most $4 k + 1$ vertices
where $k \ne 2$ contains a directed Hamiltonian circuit.
\end{conjecture}

Jackson's bound is best possible because of the disjoint union 
of two $k$-diregular tournaments on $2 k+1$ vertices.

\begin{claim}{Counterexample}

Computer search for non-hamiltonian circulant digraphs found the
following counterexample.

Let $H$ be the circulant digraph on $12$ veritces $\{0, 1, 2 \ldots 11\}$
and edges $i$ is adjacent to $\{i + 2 \mod 12, i+3 \mod 12, i+8 \mod 12\}$.
\end{claim}
\begin{proof}
By construction and inspection it is $3$-diregular and oriented.

Non-hamiltonicity was noted by results of the CAS sage \cite{sage1}
and since the digraph is small exhaustive search with \cite{sage1}
confirmed it is non-hamiltonian.

From discussion we learned $H$ appears in \cite{locke1} as
$\operatorname{Cay}(\mathbb{Z}_{12} ; 2 , 3 ,8)$ on p. 2.
Non-hamiltonicity is proved in \cite{locke1} p. 6 Theorem 4.6.
The second non-hamiltonian circulant digraph on $12$ vertices in 
\cite{locke1} is not oriented.

$H$ is on $12$ vertices, less than the bound $4 \cdot 3 + 1 = 13$
so satisfies the hypotheses for Jackson's conjecture, and not the 
conclusion. $H$ is counterexample to a weaker conjecture with
$ 4 k + 1$ replaced by $4 k$.

Computer search for circulant counterexamples didn't find
others on less than $21$ vertices.

By coincidence $\operatorname{Cay}(\mathbb{Z}_{12} ; 2 , 3 ,8)$ is
a counterexample to Adam's conjecture as shown in \cite{jira1}.
\begin{figure}
\includegraphics{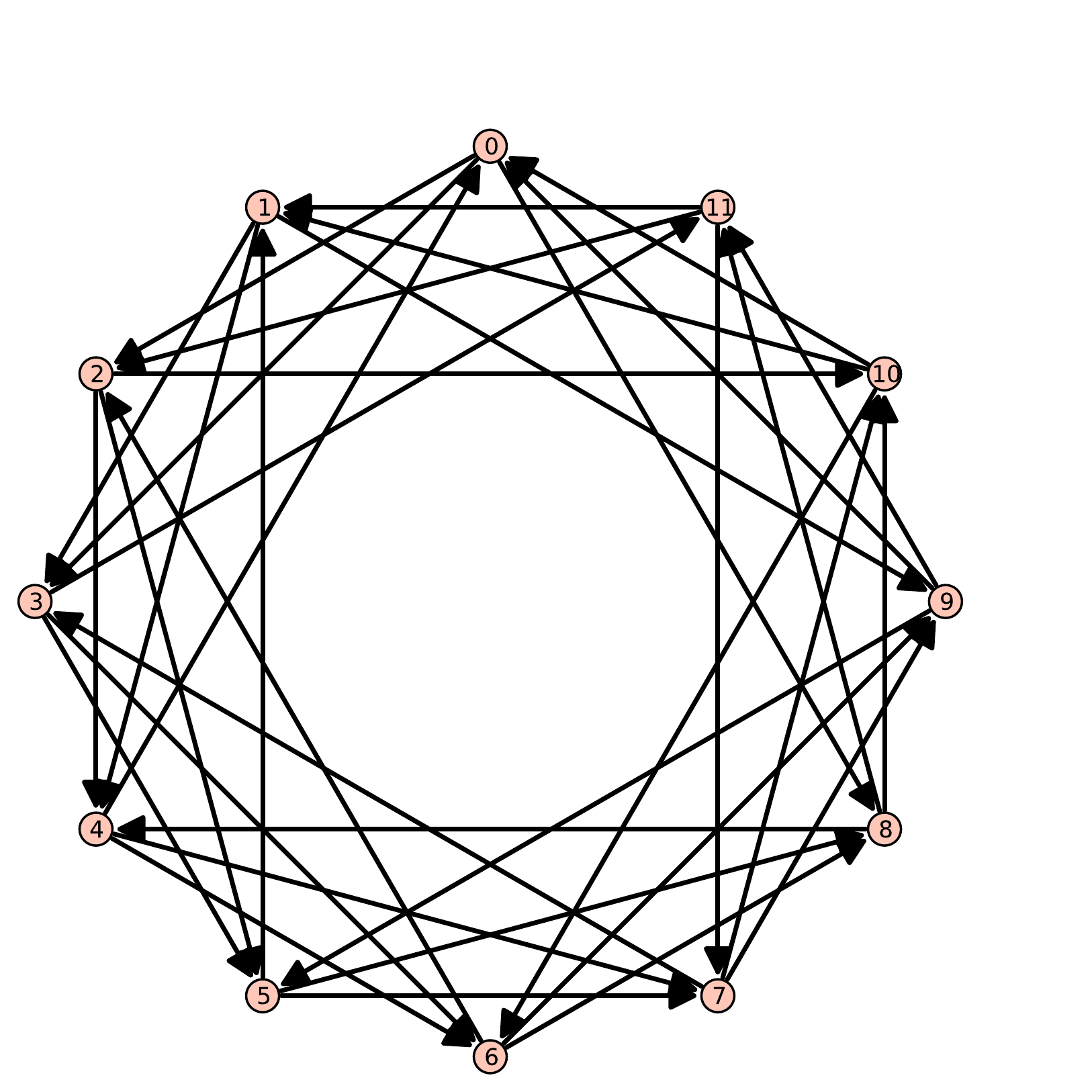}
\caption{The digraph $H$}
\label{fig:digraph}
\end{figure}
\end{proof}

\section*{Acknowledgements}
\thispagestyle{empty}

We thank Nathann Cohen for his help.

\clearpage


\end{document}